\documentclass[pdf,aos,preprint]{imsart}
\usepackage{amsmath}
\usepackage{amssymb}
\usepackage{amsthm}
\usepackage{graphicx}
\usepackage{color}
\usepackage{hyperref}
\usepackage{verbatim} 
\usepackage{arydshln}

\newtheorem{thm}{Theorem}[section] 
\newtheorem{lem}{Lemma}[section] 
\newtheorem{cor}{Corollary}[section]

\newtheorem{definition}{Definition}[section]

\newcommand{\bed}{\begin{definition}}
\newcommand{\eed}{\end{definition}}

\newcommand{\bitem}{\begin{itemize}}
\newcommand{\eitem}{\end{itemize}}

\newcommand{\beqn}{\begin{equation}}
\newcommand{\eeqn}{\end{equation}}
\newcommand{\balign}{\begin{align}}
\newcommand{\ealign}{\end{align}}

\newcommand{\cG}{{\cal G}}

\newcommand{\cL}{{\cal L}}

\newcommand{\beq}{\begin{equation}}
\newcommand{\eeq}{\end{equation}}

\allowdisplaybreaks
\begin{document} 

\begin{frontmatter}

\title{A Sharp Lower Bound for Mixed-membership Estimation}

\begin{aug}
\author{\fnms{Jiashun} \snm{Jin}\thanksref{t1}\ead[label=e1]{jiashun@stat.cmu.edu}} and 
\author{\fnms{Zheng Tracy} \snm{Ke}\thanksref{t2}\ead[label=e2]{zke@galton.uchicago.edu}}

\affiliation{Carnegie Mellon University\thanksmark{t1} and   
University of Chicago\thanksmark{t2}} 
    
         \address{J. Jin\\
		Department of Statistics\\
		Carnegie Mellon University\\
		Pittsburgh, Pennsylvania, 15213\\
		USA\\
		\printead{e1}}
	
	\address{Z. Ke\\
		Department of Statistics\\
		University of Chicago\\
		Chicago, Illinois, 60637\\
		USA\\
		\printead{e2}}
\end{aug}

\begin{abstract} 
Consider an undirected network with $n$ nodes and $K$ perceivable communities, where some nodes may have mixed memberships.    We assume that for each node $1 \leq  i \leq n$, there is a probability mass function 
$\pi_i$ defined over  $\{1, 2, \ldots, K\}$ such that 
\[
\pi_i(k) = \mbox{the weight of node $i$ on community $k$}, \qquad 1 \leq k \leq K.  
\]  
The goal is to  estimate $\{\pi_i, 1 \leq i \leq n\}$ (i.e., membership estimation).

We model the network with the {\it degree-corrected mixed membership (DCMM)} model  \cite{Mixed-SCORE}.   Since for many natural networks, the degrees have an approximate  power-law tail, we allow {\it severe degree heterogeneity} in our model.

For any membership estimation $\{\hat{\pi}_i, 1 \leq i \leq n\}$,  since each $\pi_i$ is a probability mass function, it is natural to measure the errors  by  
the average $\ell^1$-norm 
\[
\frac{1}{n} \sum_{i = 1}^n  \| \hat{\pi}_i - \pi_i\|_1.  
\]  
We also consider a  variant  of the $\ell^1$-loss, where each $\|\hat{\pi}_i - \pi_i\|_1$ is re-weighted by the degree parameter $\theta_i$ in DCMM (to be introduced).

We present a sharp lower bound.  We also show that such a lower bound is achievable under a broad situation. More discussion in this vein is continued in our forthcoming manuscript \cite{SCORE-OptimalRate}.  

The results are very different from those on community detection. For community detection, the focus is on the special case where all $\pi_i$ are degenerate; the goal is clustering, so Hamming distance is the natural choice of loss function, and the rate can be exponentially fast. 
The setting here is broader and more difficult: it is more natural to use the $\ell^1$-loss,  and the rate is only polynomially fast.    
   
\end{abstract}
%

\end{frontmatter}

\section{Introduction}

Consider an undirected network ${\cal N} = (V, E)$, where $V = \{1, 2, \ldots, n\}$ is the set of nodes and $E$ is the set of (undirected) edges. Let $A\in\mathbb{R}^{n,n}$ be the adjacency matrix where
\[
A(i,j) = \left\{ 
\begin{array}{ll}
1,  &\qquad \mbox{if nodes $i$ and $j$ have an edge},  \\
0,  &\qquad \mbox{otherwise}. \\   
\end{array} 
\right. \qquad 1 \leq i, j \leq n, 
\] 
The diagonals of $A$ are zero since we do not allow for self-edges. 
Suppose the network has $K$ perceivable communities (i.e., clusters)
\[   
{\cal C}_1,{\cal C}_2, \ldots, {\cal C}_K,
\]
and a node may belong to more than one cluster (i.e., mixed memberships). For each node $1\leq i\leq n$, suppose there exists a Probability Mass Function (PMF)  
$\pi_{i} = (\pi_{i}(1), \pi_i(2), \ldots, \pi_{i}(K))' \in \mathbb{R}^{K}$   
such that
\[
\mbox{$\pi_i(k)$ is the ``weight" of node $i$ on ${\cal C}_k$}, \qquad 1\leq k\leq K.  
\]
We call node $i$ ``pure"  if $\pi_i$ is degenerate (i.e., one entry is $1$ and the other entries are $0$)  and  ``mixed" otherwise. The   
primary interest is to estimate $\pi_i$, $1 \leq i  \leq n$. 
  
Estimating mixed memberships is a problem of great interest in social network analysis  \cite{airoldi2009mixed,Tensor,Mixed-SCORE,JiZhuMM}.
Take the Polbook network \cite{polbook} for example.  Each node is a book on US politics for sale in Amazon.com, and  there is an edge between two nodes if they are frequently co-purchased.
Jin et al. (2017) \cite{Mixed-SCORE} modeled this network with a two-community (``Conservative" and ``Liberal") mixed membership model, 
where the estimated mixed membership of a node describes 
how much weight this book puts on ``Conservative" and ``Liberal".

We are interested in the optimal rate of convergence associated with membership estimation. 
Below, we introduce a model and present a sharp lower bound. We show that the lower bound 
is achievable in a broad class of situations where we allow severe degree heterogeneity.

\subsection{Model}
Consider the degree-corrected mixed membership (DCMM) model \cite{Mixed-SCORE}. Recall that $A$ is the adjacency matrix. DCMM assumes that 
\beq \label{mod-DCMM0} 
\{A(i,j): 1\leq i<j\leq n\} \mbox{ are  independent Bernoulli variables},
\eeq
where the Bernoulli parameters are different. For a symmetric non-negative matrix $P\in\mathbb{R}^{K,K}$ and a vector $\theta=(\theta_1,\theta_2,\ldots,\theta_n)'$, where $\theta_i>0$ is the {\it degree heterogeneity parameter} of node $i$, DCMM models 
\beq \label{mod-DCMM} 
\mathbb{P}\big(A(i,j) = 1\big) =  \theta_i \theta_j \cdot \pi_i' P\pi_j, \qquad 1\leq i<j\leq n. 
\eeq 
To ensure model identifiability, we assume 
\beq  \label{mod-identifiability}
\mbox{$P$ is non-singular and have unit diagonals}. 
\eeq

We now calibrate DCMM with a matrix form. 
Introduce the two matrices $\Theta=\mathrm{diag}(\theta_1,\theta_2,\ldots,\theta_n)\in\mathbb{R}^{n,n}$ and $\Pi=[\pi_1,\pi_2,\ldots,\pi_n]'\in\mathbb{R}^{n,K}$. Then, 
\[
A = \underbrace{[\Omega - \mathrm{diag}(\Omega)]}_{\text{``signal"}} + \underbrace{W}_{\text{``noise"}}, \qquad \Omega =\Theta\Pi P\Pi'\Theta, \qquad W = A - \mathbb{E}[A]. 
\]
Here $\Omega$ is a low-rank matrix ($\mathrm{rank}(\Omega)=K$) containing Bernoulli parameters and $W$ is a generalized Wigner matrix. 

DCMM can be viewed as extending the {\it mixed membership stochastic block model} (MMSB) \cite{airoldi2009mixed}  to accommodate degree heterogeneity, and can also be viewed as extending the {\it degree-corrected block model} (DCBM) \cite{DCBM} to accommodate mixed memberships. DCMM is similar to the {\it overlapping continuous community assignment model} (OCCAM) \cite{JiZhuMM}, where the difference is that DCMM regards each membership vector $\pi_i$ as a PMF with a unit $\ell^1$-norm while OCCAM models that each $\pi_i$ has a unit $\ell^2$-norm (which seems hard to interpret).

{\bf Remark}. The identifiability condition of DCMM is different from that of DCBM. In DCBM, even when $P$ is singular, the model can still be identifiable. However, in DCMM, since there are many more free parameters, the full rank assumption of $P$ is required for identifiability.

{\bf An example}. Let's look at an example with $K=2$ and 
\[
P = \begin{pmatrix} a & b\\ b & c\end{pmatrix}. 
\]
If nodes $i$ and $j$ are both pure nodes, then there are three cases:
\[
\mathbb{P}\big(A(i,j) = 1\big) = \theta_i\theta_j \begin{cases} 
a, & i,j\in {\cal C}_1,\\
c, & i,j\in {\cal C}_2,\\
b, & i\in {\cal C}_1,j\in{\cal C}_2 \; \text{or}\; i\in {\cal C}_2,j\in{\cal C}_1.
\end{cases}
\]
As a result, in the special case with all nodes being pure, the ``signal" matrix $\Omega$ has the form
\[
\Omega =   
\begin{bmatrix}
\theta_1  &                &   &     \\
                &                & \ddots &     \\
   &   &  & \theta_n  \\ 
\end{bmatrix}
\underbrace{\begin{bmatrix}
1 & 0\\ 
0 & 1\\
\vdots & \vdots  \\ 
1 & 0 \\ 
\end{bmatrix}
\begin{bmatrix}
a  &   b     \\
b  &   c     \\ 
\end{bmatrix}
\begin{bmatrix}
1 & 0 & \cdots & 1\\ 
0 & 1 & \cdots & 0
\end{bmatrix}}_{\Pi P \Pi'}
\begin{bmatrix}
\theta_1  &                &   &     \\
                &                & \ddots &     \\
   &   &  & \theta_n  \\ 
\end{bmatrix}, 
\]
where the matrix $\Pi P\Pi'$ can be shuffled to a block-wise constant matrix by some unknown permutation: 
\[
\Pi P\Pi' = \begin{bmatrix}
a & b & a & b & a  \\
b & c & b & c & b  \\
a & b & a & b & a \\
b & c & b & c & b  \\
a & b & a & b & a  \\
\end{bmatrix}
\overset{permute}{\longrightarrow}
\left[
\begin{array}{ccc;{2pt/2pt}cc}
a & a & a & b & b    \\
a & a & a & b & b    \\
a & a & a & b & b    \\ \hdashline[2pt/2pt]
b & b & b & c & c    \\
b & b & b & c & c    \\
\end{array}\right]. 
\]
In general cases where the nodes have mixed memberships, $\Omega$ has a similar form, except that $\Pi$ is no longer a matrix of $0$'s and $1$'s and $\Pi P\Pi'$ can no longer be shuffled to a block-wise constant matrix. 
\[
\Omega =   
\begin{bmatrix}
\theta_1  &                &   &     \\
                &                & \ddots &     \\
   &   &  & \theta_n  \\ 
\end{bmatrix}
\underbrace{\begin{bmatrix}
0.8 & 0.2\\ 
0 & 1\\
\vdots & \vdots  \\ 
0.7 & 0.3 \\ 
\end{bmatrix}
\begin{bmatrix}
a  &   b     \\
b  &   c     \\ 
\end{bmatrix}
\begin{bmatrix}
0.8 & 0 & \cdots & 0.7\\ 
0.2 & 1 & \cdots & 0.3
\end{bmatrix}}_{\Pi P \Pi'}
\begin{bmatrix}
\theta_1  &                &   &     \\
                &                & \ddots &     \\
   &   &  & \theta_n  \\ 
\end{bmatrix}.
\]

\subsection{Loss functions}
Given estimators $\hat{\Pi}=[\hat{\pi}_1,\hat{\pi}_2,\ldots,\hat{\pi}_n]'$, 
since each $\pi_i$ is a PMF, it is natural to measure the (unweighted) $\ell^1$-estimation error: 
\beq \label{uwLoss}
{\cal H}(\hat{\Pi},\Pi) = n^{-1}\sum_{i=1}^n \|\hat{\pi}_i-\pi_i\|_1. 
\eeq
We also consider a variant of the $\ell^1$-error where $\|\hat{\pi}_i-\pi_i\|_1$ is reweighed by the degree parameter $\theta_i$. 
Write $\bar{\theta}=n^{-1}\sum_{i=1}^n\theta_i$, $\theta_{\max}=\max_{1\leq i\leq n}\theta_i$, and $\theta_{\min}=\min_{1\leq i\leq n}\theta_i$. Define the degree-weighted $\ell^1$-estimation error as
\beq \label{wLoss}
\cL(\hat{\Pi},\Pi) = n^{-1}\sum_{i=1}^n (\theta_i/\bar{\theta})^{1/2} \|\hat{\pi}_i-\pi_i\|_1. 
\eeq

When $\theta_{\max}/\theta_{\min}$ is bounded, the above loss functions are equivalent in the sense that for a constant $C> 1$, 
\[
C^{-1}{\cal H}(\hat{\Pi},\Pi)\leq \cL(\hat{\Pi},\Pi)\leq C{\cal H}(\hat{\Pi},\Pi). 
\]
However, when there is {\it severe degree heterogeneity} (i.e., $\theta_{\max}/\theta_{\min}\gg 1$), the weighted $\ell^1$-loss is more convenient to use: The minimax rate for ${\cal H}(\hat{\Pi},\Pi)$ depends on all $\theta_i$ in a complicated form, but the minimax rate of ${\cal L}(\hat{\Pi},\Pi)$ is a simple function of $\bar{\theta}$.

{\bf Remark}. The weights in \eqref{wLoss} are motivated by the study of an {\it oracle case} where all true parameters of DCMM are known except for $\pi_i$ of one node $i$. 
In this case, 
there exists an {\it oracle estimator} $\hat{\pi}_{i0}$ such that 
\[
\|\hat{\pi}_{i0}- \pi_i\|_1=(\theta_i/\bar{\theta})^{-1/2}\cdot O((n\bar{\theta}^2)^{-1/2})
\] with high probability. 
It motivates us to re-weight $\|\hat{\pi}_i-\pi_i\|$ by $(\theta_i/\bar{\theta})^{1/2}$. 


\subsection{Lower bound}
In the asymptotic analysis, we fix $K$ and $P\in\mathbb{R}^{K,K}$ and let $(\Theta,\Pi)$ change with $n$. 
Our results take the form: For each $\theta$ in a broad class ${\cal Q}^*_n(K,c)$ (to be introduced), we provide a minimax lower bound associated with a class of $\Pi$. 

Given $\theta\in\mathbb{R}^n_+$, let $\theta_{(1)}\leq \theta_{(2)}\leq \ldots \leq \theta_{(n)}$ be the sorted values of $\theta_i$'s. For a constant $c\in (0,1/K)$, introduce 
\beq \label{theta-class}
{\cal Q}^*_n(K,c) = \bigl\{\theta\in\mathbb{R}_+^n: \bar{\theta}\geq n^{-1/2}\log(n),\; \theta_{(cKn)}\geq n^{-1/2}\log(n)\bigr\}   
\eeq

Denote by $\cG_n(K)$ the set of  all matrices $\Pi\in\mathbb{R}^{n,K}$ such that each row $\pi_i$ is a PMF. 
Given $\Pi\in \cG_n(K)$, let 
\[
{\cal N}_k=\{1\leq i\leq n: \pi_i=e_k\}, \qquad {\cal M}=\{1,2,\ldots,n\}\setminus ({\cal N}_1\cup\ldots\cup {\cal N}_K),
\]
where $e_1,e_2,\ldots,e_K$ are the standard bases of $\mathbb{R}^K$. It is seen that ${\cal N}_k$ is the set of pure nodes of community $k$ and ${\cal M}$ is the set of mixed nodes. 
Fix $(K,c)$ and an integer $L_0\geq 1$. Introduce
\begin{align}
&\widetilde{\cG}_n(K,c, L_0;\theta)= \Bigl\{\Pi\in\cG_n(K):  \cr 
&\begin{array}{cl}
&|{\cal N}_k|\geq cn, \mbox{for } 1\leq k\leq K;\\
&\sum_{i\in{\cal N}_k}\theta^2_i\geq c\|\theta\|^2, \mbox{for } 1\leq k\leq K; \\
&\mbox{there is $L\leq L_0$, a partition ${\cal M}=\cup_{\ell=1}^{L}{\cal M}_\ell$, PMF's $\gamma_1,..., \gamma_{L}$},\\
&\mbox{where $\min_{j\neq \ell}\|\gamma_j-\gamma_\ell\|\geq c$, $\min_{1\leq \ell\leq L,1\leq k\leq K}\|\gamma_\ell-e_k\|\geq c$},\\
&\mbox{such that } |{\cal M}_\ell|\geq c|{\cal M}|\geq \frac{\log^3(n)}{\bar{\theta}^{2}}, \max_{i\in{\cal M}_\ell}\|\pi_i - \gamma_\ell\| \leq \frac{1}{\log(n)}\Bigr\}.
\end{array}
\end{align}


\begin{thm}[Lower bound of the weighted $\ell^1$-error] \label{thm:LB}
Fix $K\geq 2$, $c\in (0,1/K)$, and a nonnegative symmetric matrix $P\in\mathbb{R}^{K,K}$ that satisfies \eqref{mod-identifiability}. 
Suppose the DCMM model \eqref{mod-DCMM0}-\eqref{mod-DCMM} holds. As $n\to\infty$, there are constants $C_0>0$ and $\delta_0\in (0,1)$ such that, for any $\theta\in {\cal Q}_n^*(K,c)$, 
\[
\inf_{\hat{\Pi}}\sup_{\Pi\in\widetilde{\cG}_n(K, c, L_0;\theta)} \mathbb{P}\Bigl( \cL(\hat{\Pi}, \Pi) \geq \frac{C_0}{\sqrt{n\bar{\theta}^2}}\Bigr) \geq \delta_0. 
\]
\end{thm}

When $\theta_{\max}\leq C\theta_{\min}$, the unweighted and weighted $\ell^1$-errors are equivalent, and we also have a lower bound for the unweighted $\ell^1$-error:  
\begin{cor}[Lower bound of the unweighted $\ell^1$-error]  \label{cor:LB0}
Suppose the conditions of Theorem~\ref{thm:LB} hold. As $n\to\infty$, there are constants $C_1>0$ and $\delta_0\in (0,1)$ such that, for any $\theta\in {\cal Q}_n^*(K,c)$ satisfying $\theta_{\max}\leq C\theta_{\min}$, 
\[
\inf_{\hat{\Pi}}\sup_{\Pi\in\widetilde{\cG}_n(K, c, L_0;\theta)} \mathbb{P}\Bigl( {\cal H}(\hat{\Pi}, \Pi) \geq \frac{C_1}{\sqrt{n\bar{\theta}^2}}\Bigr) \geq \delta_0. 
\]
\end{cor}


{\bf Remark}. Our results allow for {\it severe degree heterogeneity}: for $\theta\in{\cal Q}^*_n(K,c)$, it is possible that $\theta_{\max}/\theta_{\min}\gg 1$. In addition, we allow for {\it sparse networks} because $\theta\in{\cal Q}^*_n(K,c)$ only requires that the average node degree grows with $n$ in a logarithmic rate.   

{\bf Remark}. The lower bounds here are different from those on community detection \cite{zhang2016minimax, gao2016community}. For community detection, the focus is on the special case where all $\pi_i$ are degenerate; the goal is clustering, so Hamming distance is the natural choice of loss function, and the rate can be exponentially fast. 
The setting here is broader and more difficult: it is more natural to use the $\ell^1$-loss,  and the rate is only polynomially fast.

\subsection{Achievability}
Jin et al. \cite{Mixed-SCORE} proposed a method {\it Mixed-SCORE} for estimating $\pi_i$'s. The Mixed-SCORE is a fast and easy-to-use spectral approach, and can be viewed as an extension of Jin's SCORE \cite{SCORE,JiJin,SCORE+}. However, SCORE is originally designed for community detection,  
and to extend it to membership estimation, we need several innovations; see \cite{SCORE, Mixed-SCORE} for details. 
 It turns out that Mixed-SCORE is also rate-optimal.

The following theorem follows directly form Theorem 1.2 of \cite{Mixed-SCORE}:
\begin{thm}[Upper bound] \label{thm:UB}
Fix $K\geq 2$, $c\in (0,1/K)$, and a nonnegative symmetric irreducible matrix $P\in\mathbb{R}^{K,K}$ that satisfies \eqref{mod-identifiability}. 
Suppose the DCMM model \eqref{mod-DCMM0}-\eqref{mod-DCMM} holds. Let $\hat{\Pi}$ be the Mixed-SCORE estimator. As $n\to\infty$, for any $\theta\in {\cal Q}_n^*(K,c)$ with $\theta_{\max}\leq C\theta_{\min}$ and any $\Pi\in\widetilde{\cG}_n(K, c,L_0)$, with probability $1-o(n^{-3})$, 
\[
\cL( \hat{\Pi}, \Pi)\leq C{\cal H}(\hat{\Pi}, \Pi) \leq \frac{C\log(n)}{\sqrt{n\bar{\theta}^2}}. 
\]
\end{thm}

In the case that $\theta_{\max}\leq C\theta_{\min}$, the upper bound and lower bound have matched, and the minimax rate of convergence for both weighted and unweighted $\ell^1$-errors is 
\[
(n\bar{\theta}^2)^{-1/2}, \qquad \mbox{up to a multiple-$\log(n)$ factor}. 
\]

For more general settings where $\theta_{\max}/\theta_{\min}$ is unbounded, in a forthcoming manuscript Jin and Ke \cite{SCORE-OptimalRate}, we demonstrate that
\begin{itemize}
\item The minimax rate of convergence for the weighted $\ell^1$-loss ${\cal L}(\hat{\Pi},\Pi)$ is still $(n\bar{\theta}^2)^{-1/2}$, up to a multiple-$\log(n)$ factor. 
\item The minimax rate of convergence for the unweighted $\ell^1$-loss ${\cal H}(\hat{\Pi},\Pi)$ depends on individual $\theta_i$'s in a more complicated form. 
\item Mixed-SCORE achieves the minimax rate for a broad range of settings. 
\end{itemize} 


At the heart of the upper bound arguments is  some new {\it node-wise large deviation bounds} 
we derived; see our forthcoming manuscript \cite{SCORE-OptimalRate}. 
On a high level, the technique is connected to the post-PCA entry-wise bounds in Jin {\it et al}.  \cite{JKW}  and Ke and Wang \cite{ke2017new}, but is for very different settings. The main interest of \cite{JKW} is on gene microarray analysis, where we discuss three interconnected problems:  subject clustering, signal recovery, and global testing; see also Jin and Wang \cite{IFPCA} on IF-PCA. The main interest of \cite{ke2017new} is on topic estimation in text mining.  

As far as we know,   Jin {\it et al}.   \cite{JKW} is the first paper that has carefully studied post-PCA entry-wise bounds. The bounds are crucial for obtaining sharp bounds on the  clustering errors by PCA approaches.

\section{Proof of Theorem~\ref{thm:LB}}
We introduce a subset of $\widetilde{\cG}_n(K,c,L_0;\theta)$:
\[
\begin{array}{rl}
\cG^*_n(K,c;\theta)= \bigl\{\Pi\in\cG_n(K):  &|{\cal N}_k|\geq cn, \mbox{for } 1\leq k\leq K;\\
&\sum_{i\in{\cal N}_k}\theta^2_i\geq c\|\theta\|^2, \mbox{for } 1\leq k\leq K; \\
&\|\pi_i-(1/K){\bf 1}_K\| \leq 1/\log(n), \mbox{for }i\in {\cal M}\bigr\}.
\end{array}
\]
Since $\cG_n^*(K,c;\theta)\subset \widetilde{\cG}_n(K,c,L_0;\theta)$, for any estimator $\hat{\Pi}$, 
\[
\sup_{\Pi\in\widetilde{\cG}_n(K, c, L_0;\theta)} \mathbb{P}\Bigl( \cL(\hat{\Pi}, \Pi) \geq \frac{C_0}{\sqrt{n\bar{\theta}^2}}\Bigr) \geq \sup_{\Pi\in \cG^*_n(K, c;\theta)} \mathbb{P}\Bigl( \cL(\hat{\Pi}, \Pi) \geq \frac{C_0}{\sqrt{n\bar{\theta}^2}}\Bigr). 
\]
Hence, it suffices to prove the lower bound for $\Pi\in \cG_n^*(K,c;\theta)$. 

We need the following lemma, which is adapted from Theorem 2.5 of \cite{tsybakov2009introduction}. We recall that $\cG_n(K)$ is the set of all matrices $\Pi\in\mathbb{R}^{n,K}$ each row of which is a PMF in $\mathbb{R}^K$.
\begin{lem} \label{lem:LBbasic}
For any subset ${\cal G}_n^*\subset {\cal G}_{n}(K)$, 
if there exist $\Pi^{(0)},\Pi^{(1)},\ldots,\Pi^{(J)}\in{\cal G}_n^*$ such that:
\begin{itemize}
\item[(i)] ${\cal L}(\Pi^{(j)}, \Pi^{(k)})\geq 2C_0 s_n$ for all $0\leq j\neq k\leq J$, 
\item[(ii)] $\frac{1}{J+1}\sum_{j=0}^J KL(\mathcal{P}_j,\mathcal{P}_0)\leq \beta \log(J)$,  
\end{itemize} 
where $C_0>0$, $\beta\in (0, 1/8)$, $\mathcal{P}_j$ denotes the probability measure associated with $\Pi^{(j)}$, and $KL(\cdot,\cdot)$ is the Kullback-Leibler divergence, then
\[
\inf_{\hat{\Pi}}\sup_{\Pi\in{\cal G}_{n}^*}\mathbb{P}\Bigl( \cL(\hat{\Pi},\Pi)\geq C_0s_n\Bigr) \geq \tfrac{\sqrt{J}}{1+\sqrt{J}}\Big(1-2\beta-\sqrt{\tfrac{2\beta}{\log(J)}}\Big). 
\]
As long as $J\to\infty$ as $n\to\infty$, the right hand side is lower bounded by a constant.  
\end{lem}
By Lemma~\ref{lem:LBbasic}, it suffices to find $\Pi^{(0)},\Pi^{(1)}\ldots,\Pi^{(J)}\in{\cal G}_n^*(K,c)$ that satisfy the requirement of Lemma~\ref{lem:LBbasic}. 
Below, we first consider the case $K=2$ and then generalize the proofs to $K\geq 3$. 

\subsection{The case $K=2$}
We re-parametrize the model by defining $a\in (0,1]$ and $\gamma=(\gamma_1,\ldots,\gamma_n)\in [-1,1]^n$ through
\beq \label{LBproof-1}
P = \begin{bmatrix} 1& 1-a\\1-a & 1\end{bmatrix}, \qquad \pi_i = \Big( \frac{1+\gamma_i}{2},\; \frac{1-\gamma_i}{2}\Big)', \qquad 1\leq i\leq n. 
\eeq
Since there is a one-to-one mapping between $\Pi$ and $\gamma$, we instead construct $\gamma^{(0)},\gamma^{(1)},\ldots,\gamma^{(n)}$. 
Without loss of generality, we assume $\theta_1\geq \theta_2\geq\ldots\geq\theta_n$. 
Let $n_1=\lfloor cn\rfloor$ and $n_0=n-2n_1$. Introduce
\[
\gamma^{*}= \bigl( \underbrace{0, 0, \cdots, 0}_{n_0}, \underbrace{1, 1, \cdots, 1}_{n_1}, \underbrace{-1,-1,\cdots,-1}_{n_1} \bigr)'. 
\]
Note that $\gamma^*_i\in\{\pm 1\}$ implies that node $i$ is a pure node and $\gamma^*_i=0$ indicates that $\pi_i^*=(1/2,1/2)$. 
From the Varshamov-Gilbert bound for packing numbers \cite[Lemma 2.9]{tsybakov2009introduction}, there exist $J_0\geq 2^{n_0/8}$ and $\omega_*^{(0)}, \omega_*^{(1)},\ldots,\omega_*^{(J_0)}\in \{0,1\}^{n_0}$ such that $\omega_*^{(0)}=(0,0,\ldots,0)'$ and $\|\omega_*^{(j)}-\omega_*^{(\ell)}\|_1 \geq n_0/8$, for all $0\leq j\neq \ell\leq J_0$. Let $J = 2J_0$, $\omega^{(0)}=\omega_*^{(0)}$, and $\omega^{(2\ell\pm 1)} = \pm \omega_*^{(\ell)}$ for  $1\leq \ell\leq J_0$. Then, the resulting $\omega^{(0)},\omega^{(1)},\ldots,\omega^{(J)}$ satisfy that:
\begin{itemize}
\item[(a)] $\min_{0\leq j\neq \ell\leq J}\|\omega^{(j)}-\omega^{(\ell)}\|_1\geq  n_0/8$.  
\item[(b)] For any real sequence $\{h_i\}_{i=1}^n$, $\sum_{\ell=0}^J \sum_{i=1}^{n_0}h_i\omega_i^{(\ell)}=0$. 
\end{itemize} 
For a properly small constant $c_0>0$ to be determined, letting $\delta_n=c_0(n\bar{\theta})^{-1/2}$, we construct $\gamma^{(0)}, \gamma^{(1)},\ldots,\gamma^{(J)}$ by  
\beq \label{LBproof-gamma}
\gamma^{(\ell)}=\gamma^* + \delta_n\bigl(v\circ\omega^{(\ell)},\; \underbrace{0,0,\ldots,0}_{2n_1}\bigr),
 \quad \mbox{with}\;\; v=\Bigl(\frac{1}{\sqrt{\theta_1}}, \frac{1}{\sqrt{\theta_2}},\ldots, \frac{1}{\sqrt{\theta_{n_0}}}\Bigr).  
\eeq 
We then use the one-to-one mapping \eqref{LBproof-1} to obtain $\Pi^{(0)},\Pi^{(1)},\ldots,\Pi^{(J)}$. To check that each $\Pi^{(\ell)}$ belongs to ${\cal G}^*_n(K, c;\theta)$, we notice that $\|\pi_i^{(\ell)}-(1/2,1/2)\|=O(\theta_i^{-1/2}\delta_n)=O(\theta_{n_0}^{-1/2}\delta_n)$ for $1\leq i\leq n_0$; $\theta_{n_0}$ is the $(2cn)$-smallest value of $\theta_1,\ldots,\theta_n$ and it satisfies that $\theta_{n_0}\geq n^{1/2}\log(n)$; additionally, $\bar{\theta}\geq n^{1/2}\log(n)$; it follows that $\|\pi_i^{(\ell)}-(1/2,1/2)\|=O(c_0/\log(n))$; hence, $\Pi^{(\ell)}\in \cG_n^*(K,c;\theta)$ as long as $c_0$ is appropriately small.


What remains is to show that the requirements (i)-(ii) in Lemma~\ref{lem:LBbasic} are satisfied for $s_n=(n\bar{\theta}^2)^{-1/2}$. Consider (i). Note that for any $0\leq j\neq \ell\leq J$,
\[
{\cal L}(\Pi^{(j)}, \Pi^{(\ell)})=\min_{\pm}\Bigl\{ \frac{1}{n \sqrt{\bar{\theta}}}\sum_{i=1}^n \sqrt{\theta_i} |\gamma^{(j)}_i\pm\gamma_i^{(\ell)}|\Bigr\}. 
\]
For ``$-$", the term in the brackets is at most $\delta_n\sum_{i=1}^{n_0}\theta_i^{-1/2}\leq n_0\delta_n\theta^{-1/2}_{n_0}=o(n)$; for ``$+$", this term is at least $4n_1\geq 4cn$. Therefore, the minimum is achieved at ``$-$". Furthermore, we have   
\[
{\cal L}(\Pi^{(j)}, \Pi^{(\ell)})=\frac{1}{n \sqrt{\bar{\theta}}}\sum_{i=1}^n \sqrt{\theta_i} |\gamma^{(j)}_i-\gamma_i^{(\ell)}|=  \frac{\delta_n}{n \sqrt{\bar{\theta}}}\|\omega^{(j)}-\omega^{(\ell)}\|_1 \geq \frac{n_0\delta_n}{8n\sqrt{\bar{\theta}}},  
\] 
where the last inequality is due to Property (a) of $\omega^{(0)},\ldots,\omega^{(J)}$. 
Since $\delta_n=c_0(n\bar{\theta})^{-1/2}$ and $n_0\geq (1-cK) n$, the right hand side is lower bounded by $ (c_0\epsilon_0/8)\cdot (n\bar{\theta}^2)^{-1/2}$. This proves (i). 

We now prove (ii). 
Note that  $KL(\mathcal{P}_\ell,\mathcal{P}_0)=\sum_{1\leq i<j\leq n}\Omega^{(\ell)}_{ij}\log(\Omega^{(\ell)}_{ij}/\Omega^{(0)}_{ij})$. Additionally, the parametrization \eqref{LBproof-1} yields that 
\beq \label{LBproof-2}
\Omega_{ij}=\theta_i\theta_j \bigl[(1-a/2)+(a/2)\gamma_i\gamma_j\bigr], \qquad 1\leq i\neq j\leq n. 
\eeq 
Since $\gamma_i^{(0)}=\gamma_i^{(\ell)}$ for all $i>n_0$, if both $i, j>n_0$, then $\Omega^{(\ell)}_{ij}=\Omega_{ij}^{(0)}$ and the pair $(i,j)$ has no contribution to $KL(\mathcal{P}_\ell,\mathcal{P})$. Therefore, we can write
\begin{align} \label{LBproof-(I+II)}
& \frac{1}{J+1}\sum_{\ell=0}^{J+1}KL(\mathcal{P}_\ell,\mathcal{P}_0)\cr
 =& \frac{1}{J+1}\sum_{\ell=0}^{J+1}\Bigl(\sum_{1\leq i< j\leq n_0} + \sum_{1\leq i\leq n_0, n_0<j\leq n}\Bigr) \Omega^{(\ell)}_{ij}\log(\Omega^{(\ell)}_{ij}/\Omega^{(0)}_{ij})\cr
\equiv & (I) + (II). 
\end{align}

First, consider $(I)$. From \eqref{LBproof-gamma} and \eqref{LBproof-2}, for all $1\leq i< j\leq n_0$, we have $\Omega_{ij}^{(0)}=\theta_i\theta_j(1-a/2)$ and 
\beq \label{LBproof-3}
\Omega_{ij}^{(\ell)}=\Omega^{(0)}_{ij}(1+\Delta^{(\ell)}_{ij}), \qquad \mbox{where}\quad \Delta^{(\ell)}_{ij}=\frac{a}{2-a}\frac{\delta_n^2}{\sqrt{\theta_i\theta_j}}\cdot \omega_i^{(\ell)}\omega_j^{(\ell)}.  
\eeq
Write $\Delta_{\max}=\max_{1\leq i<j\leq n_0,1\leq \ell\leq J}|\Delta^{(\ell)}_{ij}|$. Since $\theta_1\geq \ldots \geq \theta_{n_0}\gg n^{-1/2}$, we have $\Delta_{\max}= O(n^{1/2}\delta^2_n)=O((n\bar{\theta}^2)^{-1/2})= o(1)$. By Taylor expansion, $(1+t)\ln(1+t)=t+O(t^2)\leq 2|t|$ for $t$ sufficiently small. Combining the above gives
\begin{align*}
\Omega^{(\ell)}_{ij}\log(\Omega^{(\ell)}_{ij}/\Omega^{(0)}_{ij}) &= \Omega_{ij}^{(0)}(1+\Delta_{ij}^{(\ell)})\ln(1+\Delta_{ij}^{(\ell)})\cr
& \leq 2\Omega_{ij}^{(0)}|\Delta_{ij}^{(\ell)}|\cr
&\leq a\delta_n^2 \cdot \sqrt{\theta_i\theta_j}\cdot |\omega_i^{(\ell)}\omega_j^{(\ell)}|\cr
&\leq a\delta_n^2\sqrt{\theta_i\theta_j}, 
\end{align*}
where the third line is due to \eqref{LBproof-3} and the expression $\Omega_{ij}^{(0)}$. It follows that 
\beq \label{LBproof-4}
(I) \leq a\delta_n^2\sum_{1\leq i< j\leq n_0}\sqrt{\theta_i\theta_j}
\leq a\delta_n^2\cdot \Bigl(\sum_{1\leq i\leq n_0}\sqrt{\theta_i}\Bigr)^2. 
\eeq

Next, consider $(II)$. For $i\leq n_0$ and $j>n_0$, $\Omega_{ij}^{(0)}=\theta_i\theta_j(1-a/2)$ and 
\beq \label{LBproof-5}
\Omega_{ij}^{(\ell)}=\Omega^{(0)}_{ij}(1+\tilde{\Delta}^{(\ell)}_{ij}), \quad \mbox{with}\;\; \tilde{\Delta}^{(\ell)}_{ij}= \gamma^{(\ell)}_j\cdot \frac{a}{2-a}\frac{\delta_n}{\sqrt{\theta_i}}\cdot \omega_i^{(\ell)} \mbox{ and } \gamma_j^{(\ell)}\in\{\pm 1\}. 
\eeq 
Write $\tilde{\Delta}_{\max}=\max_{1\leq i\leq n_0,n_0<j\leq n, 1\leq \ell\leq J}|\tilde{\Delta}^{(\ell)}_{ij}|$. Similar to the bound for $\Delta_{\max}$, we have $\tilde{\Delta}_{\max}=O(n^{1/4}\delta_n)=O((n\bar{\theta})^{-1/4})=o(1)$. Also, by Taylor expansion, $(1+t)\ln(1+t)=t+t^2/2+O(|t|^3)\leq t+t^2$ for $t$ sufficiently small. Combining the above gives
\begin{align} \label{LBproof-6}
\Omega^{(\ell)}_{ij}\log(\Omega^{(\ell)}_{ij}/\Omega^{(0)}_{ij}) &=  \Omega_{ij}^{(0)}(1+\Delta_{ij}^{(\ell)})\ln(1+\Delta_{ij}^{(\ell)})\cr
&\leq \Omega_{ij}^{(0)}\tilde{\Delta}^{(\ell)}_{ij} + \Omega_{ij}^{(0)}(\tilde{\Delta}^{(\ell)}_{ij})^2. 
\end{align}
Motivated by \eqref{LBproof-6}, we first bound 
\begin{align*}
(II_1) &\equiv \frac{1}{J+1}\sum_{\ell=0}^J \sum_{i=1}^{n_0}\sum_{j=n_0+1}^n\Omega_{ij}^{(0)}\tilde{\Delta}^{(\ell)}_{ij}\cr
&= \frac{1}{J+1}\sum_{\ell=0}^J \sum_{i=1}^{n_0}\sum_{j=n_0+1}^n  \frac{a\delta_n}{2}\cdot \theta_j\gamma_j^{(\ell)}\cdot \sqrt{\theta_i}\omega_i^{(\ell)}\cr
& = \frac{a\delta_n}{2(J+1)}\Bigl(\sum_{j=n_0+1}^n\theta_j\gamma_j^{(\ell)}\Bigr)\cdot  \sum_{\ell=0}^J \sum_{i=1}^n\sqrt{\theta_i}\omega_i^{(\ell)} \cr
&=0,  
\end{align*}
where we have used Property (b) of $\omega^{(0)},\omega^{(1)},\cdots,\omega^{(J)}$. We then bound
\begin{align*}
(II_2) & \equiv \frac{1}{J+1}\sum_{\ell=0}^J \sum_{i=1}^{n_0}\sum_{j=n_0+1}^n\Omega_{ij}^{(0)}(\tilde{\Delta}^{(\ell)}_{ij})^2\cr
& = \frac{1}{J+1}\sum_{\ell=0}^J \sum_{i=1}^{n_0}\sum_{j=n_0+1}^n \frac{a^2\delta^2_n}{4-2a}\cdot \theta_j \cdot |\omega_i^{(\ell)}|\cr
&\leq \frac{a^2\delta^2_n}{4-2a} \Bigl(\sum_{j=n_0+1}^{n}\theta_j\Bigr) \cdot \max_{0\leq \ell\leq J}\|\omega^{(\ell)}\|_1\cr
&\leq \frac{a^2\delta^2_n}{4-2a}  \Bigl(\sum_{j=n_0+1}^{n}\theta_i\Bigr)\cdot n_0. 
\end{align*}
Combining the above gives
\beq \label{LBproof-7}
(II)\leq \frac{a^2\delta^2_n}{4-2a} \cdot n_0 \Bigl(\sum_{j=n_0+1}^{n}\theta_i\Bigr).
\eeq

Last, we combine \eqref{LBproof-4} and \eqref{LBproof-7}. Using the Cauchy-Schwartz inequality, $(\sum_{1\leq i\leq n_0}\sqrt{\theta_i})^2\leq n_0\sum_{i=1}^{n_0}\theta_i$. Hence, the right hand side of \eqref{LBproof-4} is upper bounded by $a\delta_n^2 \cdot n_0(\sum_{i=1}^{n_0}\theta_i)$. Furthermore, since $a\in (0,1]$, the right hand side of \eqref{LBproof-7} is upper bounded by $a\delta_n^2 \cdot n_0(\sum_{i=n_0+1}^{n}\theta_i)$. As a result, 
\[
\frac{1}{J+1}\sum_{j=0}^J KL({\cal P}_j, {\cal P}_0)\leq a\delta_n^2\cdot n_0\bigl(\sum_{i=1}^n\theta_i\bigr) = ac_0n_0,  
\]
where we have plugged in $\delta_n=c_0(n\bar{\theta})^{-1/2}$. At the same time, $\log(J)\geq [\log(2)/8]n_0$. Hence, the requirement (ii) is satisfied as long as $c_0$ is chosen appropriately small. The proof for $K=2$ is now complete.

\subsection{The case of $K\geq 3$}
The key step is to generalize the construction of $\Pi^{(0)},\Pi^{(1)},\ldots,\Pi^{(J)}$ for $K=2$ to a general $K$. 
Write $\check{P}={\bf 1}_K{\bf 1}_K'-P$. Let $\eta\in\mathbb{R}^K$ be a nonzero vector such that 
\beq \label{LBproof-8}
\eta'{\bf 1}_K=0, \qquad \eta'\check{P}{\bf 1}_K=0. 
\eeq
Such an $\eta$ always exists. 
We assume $\theta_1\geq \theta_2\geq \ldots\geq \theta_n$ without loss of generality. 
Let $n_1=\lfloor cn\rfloor$ and $n_0=n-Kn_1$. Denote by $e_1,\ldots,e_K$ the standard basis vectors of $\mathbb{R}^K$. 
Introduce
\[
\Pi^{*}= \Bigl( \underbrace{\tfrac{1}{K}{\bf 1}_K,\; \cdots,\; \tfrac{1}{K}{\bf 1}_K}_{n_0},\; \underbrace{e_1, \cdots, e_1}_{n_1},\; \cdots,\; \underbrace{e_K,\cdots,e_K}_{n_1} \Bigr)'. 
\]
Let $\omega^{(0)},\omega^{(1)},\ldots,\omega^{(J)}\in\{0,1\}^{n_0}$ be the same as above. Let $\delta_n=c_0(n\bar{\theta})^{-1/2}$ for a constant $c_0$ to be determined. Write $\Pi^*=[\pi_1^*,\ldots,\pi_n^*]'$. For each $0\leq \ell\leq J$, we construct $\Pi^{(\ell)}=[\pi_1^{(\ell)},\ldots,\pi_n^{(\ell)}]'$ by 
\[
\pi_i^{(\ell)} = \pi_i^* + \begin{cases} 
\omega_i^{(\ell)}\cdot (\delta_n /\sqrt{\theta_i}) \cdot \eta, & 1\leq i\leq n_0\\
{\bf 0}_K, & n_0+1\leq i\leq n. 
\end{cases}
\]
Same as before, we show the requirements (i)-(ii) in Lemma~\ref{lem:LBbasic} are satisfied. Note that
\[
{\cal L}(\Pi^{(j)}, \Pi^{(\ell)})= \frac{\delta_n\|\eta\|_1}{n \sqrt{\bar{\theta}}}  \|\omega^{(j)}-\omega^{(\ell)}\|_1 \geq \frac{n_0\delta_n\|\eta\|_1}{8n\sqrt{\bar{\theta}}} = (c_0\epsilon_0\|\eta\|_1/8)\cdot (n\bar{\theta}^2)^{-1/2}. 
\] 
Hence, (i) holds for $s_n=(n\bar{\theta}^2)^{-1/2}$. 

It remains is to prove (ii). Let $(I)$ and $(II)$ be defined in the same way as in \eqref{LBproof-(I+II)}. We aim to find  expressions similar to those in \eqref{LBproof-3} and \eqref{LBproof-5} and then generalize the bounds of $(I)$ and $(II)$ for $K=2$ to a general $K$.  For preparation, we first derive an expression for $\Omega_{ij}$. Introduce $\check{\pi}_i=\pi_i - \frac{1}{K}{\bf 1}_K\in\mathbb{R}^{K}$. Note that $\pi_i'{\bf 1}_K=1$. By direct calculations,  
\begin{align}  \label{LBproof-9}
\Omega_{ij} &= \theta_i\theta_j\pi_i'P\pi_j = \theta_i\theta_j(1 - \pi_i'\check{P}\pi_j)\cr
& = \theta_i\theta_j - \theta_i\theta_j (\check{\pi}_i + \tfrac{1}{K}{\bf 1}_K)'\check{P}(\check{\pi}_j + \tfrac{1}{K}{\bf 1}_K)\cr
& = \theta_i\theta_j(1- \tfrac{1}{K^2}{\bf 1}_K'\check{P}{\bf 1}_K) - \theta_i\theta_j \tfrac{1}{K}(\check{\pi}_i'\check{P}{\bf 1}_K + \check{\pi}_j'\check{P}{\bf 1}_K) - \theta_i\theta_j \check{\pi}_i'\check{P}\check{\pi}_j. 
\end{align}
 
Consider $(I)$. Since $\check{\pi}_i^{(0)}$ is a zero vector for all $1\leq i\leq n_0$, we have 
\[
\Omega_{ij}^{(0)}=\theta_i\theta_j(1-\check{a}), \quad \mbox{with}\quad \check{a}\equiv \tfrac{1}{K^2}{\bf 1}_K'\check{P}{\bf 1}_K, \qquad \mbox{if }1\leq i\neq j\leq n_0. 
\]
Furthermore, 
for $1\leq i\leq n_0$, $\check{\pi}_i^{(\ell)}\propto \eta$, where it follows from \eqref{LBproof-8} that $\eta'\check{P}{\bf 1}_K=0$. Hence, the middle two terms in \eqref{LBproof-9} are zero. As a result,  
\[
\Omega_{ij}^{(\ell)} = \theta_i\theta_j(1-\check{a}) - \frac{\delta_n^2}{\sqrt{\theta_i\theta_j}}\cdot\eta'\check{P}\eta\cdot \omega_i^{(\ell)}\omega_j^{(\ell)}, \qquad 1\leq i\neq j\leq n_0. 
\]
Combining the above, we find that for all $1\leq i\neq j\leq n_0$, 
\beq \label{LBproof-10}
\Omega_{ij}^{(\ell)}=\Omega^{(0)}_{ij}(1+\Delta^{(\ell)}_{ij}), \qquad \mbox{where}\quad \Delta^{(\ell)}_{ij}=\frac{-(\eta'\check{P}\eta)\check{a}}{1-\check{a}}\cdot \frac{\delta_n^2}{\sqrt{\theta_i\theta_j}}\cdot \omega_i^{(\ell)}\omega_j^{(\ell)}. 
\eeq
This provides a counterpart of \eqref{LBproof-3} for a general $K$. Same as before, we have the bound: $\Delta_{\max}\equiv \max_{1\leq i\neq j\leq n_0}\max_{0\leq \ell\leq J}|\Delta_{ij}^{(\ell)}|=O(n^{1/2}\delta_n^2)=o(1)$. Following the proof of \eqref{LBproof-4}, we find that
\beq \label{LBproof-11}
(I)\leq C_1 \delta_n^2\Bigl(\sum_{i=1}^{n_0}\sqrt{\theta_i}\Bigr)^2 \leq C_1n_0 \delta_n^2\Bigl(\sum_{i=1}^{n_0}\theta_i\Bigr), \qquad \mbox{where}\quad C_1= |\eta'\check{P}\eta|\cdot |\check{a}|.  
\eeq
Consider $(II)$. In this case, we need to calculate $\Omega_{ij}$, $1\leq i\leq n_0$, $n_0<j\leq n$. Recall that $\check{\pi}^{(0)}_i$ is a zero vector. Write $\{n_0+1,n_0+2,\ldots,n\}=\cup_{k=1}^K {\cal N}_k$, where ${\cal N}_k=\{1\leq j\leq n: \pi_j^{(0)}=e_k\}$. For $j\in {\cal N}_k$, it holds that $\check{\pi}_j^{(0)}+\tfrac{1}{K}{\bf 1}_K=e_k$. Combining the above with \eqref{LBproof-9}, we find that for $1\leq k\leq K$, 
\[
\Omega_{ij}^{(0)}=\theta_i\theta_j(1-\check{b}_k), \quad \mbox{with}\quad \check{b}_k \equiv \tfrac{1}{K}e_k'\check{P}{\bf 1}_K, \qquad \mbox{if }1\leq i\leq n_0, j\in {\cal N}_k.  
\]
Additionally, we have $(\check{\pi}_i^{(\ell)})'\check{P}{\bf 1}_K\propto \eta'\check{P}{\bf 1}_K=0$ and $\check{\pi}_j^{(\ell)}+\tfrac{1}{K}{\bf 1}_K=e_k$. It follows from \eqref{LBproof-9} that
\begin{align*}
\Omega_{ij}^{(\ell)} &= \theta_i\theta_j(1-\check{b}_k) - \theta_i\theta_j(\check{\pi}_i^{(\ell)})'\check{P}\check{\pi}_j^{(\ell)}\cr
& = \theta_i\theta_j(1-\check{b}_k) - \theta_i\theta_j\cdot \frac{\omega_i^{(\ell)}\delta_n}{\sqrt{\theta_i}}\cdot \eta'\check{P}(e_k - \tfrac{1}{K}{\bf 1}_K)\cr
& = \theta_i\theta_j(1-\check{b}_k) -  \theta_i\theta_j \cdot \frac{\omega_i^{(\ell)}\delta_n}{\sqrt{\theta_i}}\cdot (\eta'\check{P}e_k),
\end{align*}
where the last equality is because of $\eta'\check{P}{\bf 1}_K=0$. As a result, for $1\leq i\leq n_0$ and $j\in {\cal N}_k$, 
\beq \label{LBproof-12}
\Omega_{ij}^{(\ell)}=\Omega^{(0)}_{ij}(1+\tilde{\Delta}^{(\ell)}_{ij}), \qquad \mbox{where}\quad \tilde{\Delta}^{(\ell)}_{ij}= -(\eta'\check{P}e_k)\cdot \frac{1}{1-\check{b}_k}\frac{\delta_n}{\sqrt{\theta_i}}\cdot \omega_i^{(\ell)}. 
\eeq 
This provides a counterpart of \eqref{LBproof-5} for a general $K$. Same as before, let $\tilde{\Delta}_{\max}\equiv \max_{1\leq i\neq j\leq n_0}\max_{0\leq \ell\leq J}|\tilde{\Delta}_{ij}^{(\ell)}|$, and it is seen that $\tilde{\Delta}_{\max}=O(n^{1/4}\delta_n)=o(1)$. Again, by Taylor expansion, we have \eqref{LBproof-6}. It follows that 
\[
(II)\leq (II_1)+(II_2),
\]
where $(II_1)$ and $(II_2)$ are defined the same as before. Using \eqref{LBproof-12}, we have  
\begin{align*}
(II_1) &= \frac{1}{J+1}\sum_{\ell=0}^J \sum_{i=1}^{n_0}\sum_{j\in{\cal N}_1\cup\ldots\cup{\cal N}_K} \Omega_{ij}^{(0)}\tilde{\Delta}^{(\ell)}_{ij}\cr
&=  \frac{\delta_n}{J+1}\Bigl(\sum_{k=1}^K\sum_{j\in{\cal N}_k} -(\eta'\check{P}e_k)\theta_j\Bigr)\cdot  \sum_{\ell=0}^J \sum_{i=1}^n\sqrt{\theta_i}\omega_i^{(\ell)}\cr
&=0,
\end{align*}
where the last equality is due to Property (b) for $\omega^{(0)},\omega^{(1)},\ldots,\omega^{(J)}$. Using \eqref{LBproof-12} again, we have  
\begin{align*}
(II_2) & = \frac{1}{J+1}\sum_{\ell=0}^J \sum_{i=1}^{n_0}\sum_{j\in{\cal N}_1\cup\ldots\cup{\cal N}_K}\Omega_{ij}^{(0)}(\tilde{\Delta}^{(\ell)}_{ij})^2\cr
&\leq \frac{\delta_n^2}{1-\check{b}} \cdot  \Bigl(\sum_{1\leq k\leq K}(\eta'\check{P}e_k)^2\sum_{j\in{\cal N}_k}\theta_j\Bigr) \cdot \max_{0\leq \ell\leq J}\|\omega^{(\ell)}\|_1\cr
&\leq \frac{\delta_n^2}{1-\check{b}} \cdot \max_{1\leq k\leq K}(\eta'\check{P}e_k)^2\cdot \Bigl(\sum_{j=n_0+1}^{n}\theta_i\Bigr)\cdot n_0. 
\end{align*}
Combining the above gives
\beq \label{LBproof-13}
(II)\leq C_2n_0 \delta_n^2\Bigl(\sum_{i=n_0+1}^{n}\theta_i\Bigr), \qquad \mbox{where}\quad C_2= \frac{1}{1-\check{b}} \max_{1\leq k\leq K}(\eta'\check{P}e_k)^2.  
\eeq
We note that \eqref{LBproof-11} and \eqref{LBproof-13} server as the counterpart of \eqref{LBproof-5} and \eqref{LBproof-7}, respectively. Similarly as in the case of $K=2$, we obtain (ii) immediately. The proof for a general $K$ is now complete. 
\qed

\bibliographystyle{imsart-number}
\bibliography{network}
\end{document}